\documentstyle[12pt, amstex, amssymb]{amsart}

\def \cal #1{{\fam2 #1}}

\begin{document}
\title{\bf Support varieties of non-restricted modules
over Lie algebras of reductive groups: corrigenda and addenda}
\author{\sc
Alexander Premet}
\address{Department of Mathematics,
University of Manchester,
Oxford Road,
M13 9PL, UK}
\email{sashap@@ma.man.ac.uk}
\thanks{}
\date{}
\maketitle

J.~C.~Jantzen informed me that the proof of Lemma~3.2 in my paper
[5] is not correct in the case where $R$ is of type ${\bf C}_l$. Indeed,
the formula for
$e^{-}(t)$ at the bottom line of P.~242 should read
$$
e^{-}(t) \, =\, e_{\tilde{\alpha}}+tN_{\gamma, \tilde{\alpha}}e_{-\beta}
+\frac{1}{2}N_{\gamma, \tilde{\alpha}}N_{\gamma,-\beta}t^{2}
e_{-2\epsilon_2}.
$$
Then $\langle e^{+}(t),\, e^{-}(t)\rangle\, =\, -t^4$, and
the coefficient at $t^4$ in the expression for $b_{f}(e^{+}(t),
e^{-}(t))$ becomes $4a_{0}^{2}-4=0$ (instead of $4a_{0}+4=8$ as stated
in [5, P.~243]). So one cannot conclude that
$b_{f}(e^{+}(t), e^{-}(t))$ is a nonzero polynomial for an arbitrary
nonzero $f\in {\frak g}^*$.
In fact, if $R$ is of type ${\bf C}_l$, then
$b_{f}$ vanishes on $\cal E\times \cal E$ for some nonzero
$f\in {\frak g}^*$ (this was pointed out by H.~Kraft and N.~Wallach).

The only place where I use Lemma~3.2 is the proof of
Proposition~3.3. The purpose of this note is to derive slightly weaker
versions of Lemma~3.2 and Proposition~3.3, fill the gap in the proof of
the main result of [5] (Theorem~1.1), and
eliminate the assumption on $p$ in the formulation of
Theorem~1.1. We adopt the notation of [5].

\medskip
{\bf Lemma~3.2.} {\it Suppose that} $G\not\cong \mbox {SL}\, (2)$ {\it and}
$f\in {\frak g}^{*}\setminus \{0\}$. {\it Then}
$b_{f}|_{\, \cal E\times \cal E}
\neq 0$ {\it unless} $G\cong \mbox{Sp}\, (2l, K)$ {\it and}
$f=\langle e, \, \cdot\, \rangle$ {\it for some} $e\in \cal E$.

\begin{pf}
If $R\not\cong {\bf C}_l$ one argues as in [5].
Now assume that $p\ne 2$ and $R\cong {\bf C}_l$ where $l\ge 2$.
Let $f\in {\frak g}^{*}\setminus \{0\}$ be such that $b_{f}|_{\, \cal E\times
\cal E}\, =\, 0$.
It is well known that ${\frak g}\, =\, {\frak{sp}}\, (2l, K)$ is a simple Lie
algebra (see, e.g., [2]). Hence the $(\mbox{Ad}\, G)$--invariant bilinear
form $\langle\, \cdot
\, , \,\cdot \rangle$ is nondegenerate. It follows that $f\, =\, \langle u ,\,
\cdot \, \rangle$ for some $u\in \frak g$. Let ${\frak b}\, =\, {\frak t}
\oplus {\frak n}_+$ where ${\frak t}\, =\, \mbox{Lie} (T)$ and
$$
{\frak n}_{+}\, =\, \bigoplus_{\alpha\in R_{+}}\, K\, e_{\alpha}.
$$
By [1, P.~355], ${\frak g}\, =\, (\mbox{Ad}\, G)\cdot \frak b$. So no
generality is lost by assuming that $u\in \frak b$, that is $u=h+n$ where
$h\in \frak t$ and $n\in {\frak n}_+$.

Suppose that $h\ne 0$. As $p\ne 2$ the linear space $\frak t$ is
spanned by $h_{2\epsilon_i}$ where $1\le i\le l$. Therefore,
there exists a long root $\delta\in R_+$ such that
$\delta (h)\ne 0$. Put $v=[u, e_{\delta}]$. Then $v \ne 0$ and $f(e_{\delta})=
\langle h+n, e_{\delta} \rangle \in \langle {\frak b}, {\frak n}_{+} \rangle
=0$. As $b_{f}(\cal E, \cal E)=0$ one has
$f([e_{\delta}, e])^{2}+4 f(e_{\delta})f(e)\langle e_{\delta}, e \rangle = 0$
for any $e\in \cal E$. This gives
$$
\langle v, e \rangle = \langle u, [e_{\delta}, e] \rangle = 0
$$
whenever $e\in \cal E$. As $\cal E$ spans $\frak g$ ([5, Lemma~2.3(i)])
we get $\langle v, \frak g \rangle =0$ which contradicts the simplicity of
$\frak g$. So $h=0$ and $u\in {\frak n}_+$.

Thus $u$ is a nilpotent element of ${\frak g}$. As $p\ne 2$, $u$ has at least
one Dynkin torus $\lambda: {\bf G}_{m}\rightarrow G$
(see [4, Definition~2.4 and Theorem~2.5] and [7, Chapter~IV, (2.23)]).
Interchanging $u$ by its $(\mbox{Ad}\, G)$--conjugate if necessary we may
assume that $\lambda \subset T$. The torus $\lambda$ gives
$\frak g$ a ${\Bbb Z}$--graded Lie algebra structure
$$
{\frak g}\, =\, \bigoplus_{i\in {\Bbb Z}}\, {\frak g} (i),
$$
where ${\frak g} (i)$ is the subspace consisting of all $x\in \frak g$ with
$(\mbox{Ad}\, \lambda (t))\cdot x=t^{i} x$ for all $t\in K^*$. As $\lambda$
is a Dynkin torus for $u$ one has $u\in {\frak g} (2)$ and $\mbox{Ker}\,
\mbox{ad}\, u\subset \bigoplus_{i\ge 0}\, {\frak g} (i)$.

Let $r = \max \, \{i\in {\Bbb Z}\, |\, {\frak g} (i) \ne (0)\}$.
It is well known that there exists a parabolic subgroup $P_{-}$ of $G$ with
Levi
decomposition $P_{-}=L\cdot U_{-}$ such that $\mbox{Lie}\, (U_{-})\, =\,
\bigoplus_{i<0}\, {\frak g} (i)$ and $\mbox{Lie}\, (L)\, =\, {\frak g} (0)$.
Let $U_{0}$ be a maximal unipotent subgroup of $L$ and ${\frak u}_{0}\, =\,
\mbox{Lie}\, (U_{0})$. Then $U_{0}\cdot U_{-}$ is a maximal unipotent
subgroup of $G$, and its Lie algebra $\frak u$ equals ${\frak u}_{0}\oplus
\bigoplus_{i<0}\, {\frak g} (i)$. Obviously, ${\frak g} (r)$ is an ideal of the
Lie algebra $\frak u$. As $\frak u$ is nilpotent, ${\frak g} (r)$ intersects
with the centre of $\frak u$. As $p\ne 2$, standard properties of the root
system of type ${\bf C}_l$ ensure that the centre of $\frak u$ has
dimension $1$ and is spanned by a long root element. It follows that the set
$\cal E_{-r}\, :=\, {\frak g} (-r)\cap \cal E$ is nonzero. As $\frak g$
is a simple Lie algebra ${\frak g} (-r)$ is an irreducible $(\mbox{ad}\,
{\frak g} (0))$--module. But then
the $(\mbox{Ad}\, L)$--module ${\frak g} (-r)$ is irreducible whence
$\cal E_{-r}$ spans ${\frak g} (-r)$. As
$(\mbox{Ker}\, \mbox{ad}\, u)\cap{\frak g} (-r)\, =\, (0)$ there is
$e_{-r}\in\cal E_{-r}$ such that $[u, e_{-r}]\ne 0$.

Suppose that $r>2$. In this case $f({\frak g} (-r)) =
\langle u, {\frak g} (-r) \rangle = 0$,
and we derive that
$$
b_{f}(e_{-r}, e)=\langle u, [e_{-r}, e] \rangle^{2} =0
$$
for any $e\in \cal E$. In view of [5, Lemma~2.3(i)]
this yields $\langle [u, e_{-r}], {\frak g} \rangle = 0$
contradicting the simplicity of $\frak g$.

Thus $r=2$. Let ${\frak g} (-r)^{\prime}$ denote the subspace of all elements
in ${\frak g} (-r)$ orthogonal to $u$ with respect to $\langle\, \cdot\, ,\,
\cdot \, \rangle$.
As $\cal E$ is a Zariski closed, conical subset of $\frak g$
(see, e.g., [5, Lemma~2.1]), so is $\cal E_{-r}$. Let $d=\dim \cal E_{-r}$.
As ${\frak g} (-r)^\prime$ is a hyperplane in ${\frak g} (-r)$ one has
$\dim \cal E_{-r}\cap {\frak g (-r)}^{\prime}\ge d-1$. If $d>1$, standard
algebraic geometry says that $\cal E_{-r}\cap {\frak g} (-r)^{\prime}
\, \ne \, \{0\}$. In other words, there exists $e^{\prime}_{-r}\in
\cal E_{-r}\setminus \{0\}$ for which $\langle e^{\prime}_{-r}, u \rangle =0$.
But then arguing as before shows that $[u, e^{\prime}_{-r}]\ne 0$ and
$\langle [u, e^{\prime}_{-r}], \frak g \rangle = 0$. Since this contradicts
the simplicity of $\frak g$ we must have $d=1$.

Since $\cal E_{-r}$ is conical each irreducible component of $\cal E_{-r}$ is
a line invariant under the adjoint action of the connected group $L$. So the
irreducibility of the $(\mbox{Ad}\, L)$--module ${\frak g} (-r)$ implies that
$\dim {\frak g} (-r)=1$. As $\langle\, \cdot\,,\,\cdot\,\rangle$ defines a
non--degenerate pairing between ${\frak g} (r)$ and ${\frak g} (-r)$ the subspace
${\frak g} (r)$ is one-dimensional as well. From this it is immediate that
$u\in {\frak g} (r)$ belongs to the centre of the Lie algebra of a maximal
unipotent subgroup of $G$. As $p\ne 2$ this yields $u\in \cal E$ as required.
\end{pf}

\medskip
{\bf Proposition~3.3.}
{\it Suppose} $G$ {\it is simple and} $G\not\cong SL(2)$.  {\it Let}
$f \in {\frak g}^*$. {\it If} $G\cong \mbox{Sp}\, (2l, K)$ {\it suppose
that} $f$ {\it is not of the form} $\langle u,\, \cdot\, \rangle$
{\it with} $u\in \cal E$.
{\it Let} $x$ {\it be a nilpotent element of} $\frak g$
{\it such that} $f([x, {\frak g}])\ne 0$. {\it Then there exists} $e \in \cal E$
{\it for which} $f(e) = 0$ {\it and} $f([x, e]) \ne 0$.

\begin{pf}
One argues as in [5, PP.~243, 244] using the weaker version of [5, Lemma~3.2]
proved above. The use of that lemma is justified by the present assumption
on $f$.
\end{pf}

Let us now proceed to the proof of the main result of [5] (Theorem~1.1).
Let $\xi \in {\frak g}^*$ and let $M$ be an irreducible
$u({\frak g}, \xi)$--module. Let $z\in \cal N_{p}(\frak g)$ be such that
$\xi ([z, {\frak g}])\ne 0$. Suppose that $z\in \cal V_{\frak g}(M)$.

If $R$ is not of type ${\bf C}_l$, where $l\ge 2$, the argument in [5,
PP.~244--246] goes through. From now on assume that
$p\ne 2$ and $R\cong {\bf C}_l$ where $l\ge 2$. Then the bilinear form
$\langle\,\cdot\,,\,\cdot\,\rangle$ is nondegenerate so that $\xi=
\langle w,\cdot\,\rangle$ for some nonzero $w\in\frak g$. If $w\not\in\cal E$,
then $b_{\xi}(\cal E,\cal E)\ne 0$ (Lemma~3.2). So one argues as in [5,
PP.~244--246]. Thus in order to complete the proof it remains to treat
the case where $w\in \cal E$. No generality is lost by assuming
that $w = e_{\tilde{\alpha}}$ (recall that $\tilde{\alpha}=2\epsilon_1$).

There exists a one-dimensional torus $h_{1}(t)\subset T$
such that
$(\mbox{Ad}\, h_{1}(t))\cdot
e_{\delta}=t^{\langle \delta , \tilde{\alpha}\rangle}\cdot e_\delta$
for every $\delta\in R$ and every $t\in K^*$. Decompose $\frak g$ into weight
spaces relative to $\mbox{Ad}\, h_{1}(t)$ giving a $\Bbb Z$--grading
$$
 {\frak g}\, = \, {\frak g}(-2)\oplus {\frak g}(-1)\oplus {\frak g}(0)
\oplus {\frak g}(1)\oplus {\frak g}(2).
$$
It is easy to see that
${\frak g}(\pm 2) = K\, e_{\pm \tilde{\alpha}}$.

Let $z=\sum_{i} z_{i}$ where $z_{i}\in {\frak g} (i)$.
As $\xi ([z, {\frak g}])\ne 0$
one has $[e_{\tilde{\alpha}}, z]\ne 0$. Suppose $z_{-2}=z_{-1}=0$. Then
$z$ belongs to the parabolic subalgebra $\bigoplus_{i\ge 0}\, {\frak g} (i)$.
Since $z$ is a nilpotent element of $\frak g$ so is $z_{0}$ (this follows
from Jacobson's identity [5, (1.1)]). But then $\mbox{ad}\, z_{0}$ acts
trivially on the one-dimensional subspace $\frak g (2)$ yielding
$[z, e_{\tilde{\alpha}}]=0$. This contradiction shows that $z_{-1}\ne 0$
or $z_{-2}\ne 0$.

Suppose that $z_{-2}\ne 0$. Then $\xi (z)=\langle e_{\tilde{\alpha}},
z_{-2}\rangle \ne 0$. As $z^{[p]}=0$ the only eigenvalue of $z$ on $M$
equals $\xi (z)\ne 0$.
Let $\lambda_{z}$ be a Dynkin torus for $z$ (see [4, Definition~2.4 and
Theorem~2.5]). The subspace $\mbox{Lie}\, (\lambda_{z}) \subset \frak g$
is spanned by a toral element $H$ satisfying $[H, z]=2z$
(recall that an element $t\in \frak g$ is called {\it toral} if
$t^{[p]}=t$). Put
${\frak a}\, =\, K\, H\oplus K\, z$. Clearly, $\frak a$ is a two-dimensional
subalgebra of $\frak g$ and $K\, z$ is an ideal of $\frak a$.
But then [3, Lemma~2.9] shows that the $u(z,\xi)$--module $M$ is free
(this lemma applies as $H$ is a toral element, see [3, P.~109] for
more detail).
Since in this case $z\not\in\cal V_{\frak g}(M)$ ([5, (3.1)]) we must have
$z_{-2}=0$ and $z_{-1}\ne 0$.

Let $P_{1}=L_{1}\cdot U_{1}$ be the parabolic subgroup of $G$
such that $\mbox{Lie}\, (L_{1})\, =\, \frak g (0)$ and
$\mbox{Lie}\,(U_{1})\, =\,\bigoplus_{i>0}\, {\frak g} (i)$.
Let $L^{\prime}_1$ denote the derived subgroup of $L_1$ and $P_{1}^{\prime}
\, =\, L^{\prime}_{1}\cdot U_1$.
It is well known (and easy to check) that $L^{\prime}\cong \mbox{Sp}\,
(2l-2, K)$ and the $(\mbox{Ad}\, L^{\prime}_{1})$--module
${\frak g} (-1)$ is isomorphic to the standard $\mbox{Sp}\, (2l-2, K)$--module
of dimension $2l-2$. This means that all nonzero elements of ${\frak g} (-1)$
form a single $(\mbox{Ad}\, L^{\prime}_{1})$--orbit. In other words,
no generality is lost by assuming that $z_{-1}=e_{\gamma}$
where $\gamma=-(\epsilon_{1}+\epsilon_{2})$. Given an element
$y = t + \sum_{\alpha \in R}\,
\mu_{\alpha} e_{\alpha}$ in $\frak g$, where $\mu_{\alpha}\in K$, denote
$$
\mbox{Supp}\, (y)\, := \, \{\alpha\in R\, |\, \mu_{\alpha} \ne 0\}.
$$

The group $P^{\prime}_1$ acts trivially on the one-dimensional subspace
${\frak g} (2)$. It follows that $(\mbox{Ad}^{*} g)\cdot \xi=\xi$ for each
$g\in P^{\prime}_1$. Therefore, $P^{\prime}_1$ preserves the ideal
$I_{\xi}$ of $U(\frak g)$ whence acts by automorphisms on the quotient
algebra $u({\frak g}, \xi)=U({\frak g})/I_{\xi}$. Let $\tilde {M}$ denote the
socle of the left regular $u({\frak g} ,\xi)$--module $u({\frak g}, \xi)$, and
let $\cal V = \cal V_{\frak g}(\tilde{M})$. As $P^{\prime}_1$
preserves $\tilde{M}$ the variety $\cal V$ is
$(\mbox{Ad}\, P^{\prime}_{1})$--stable (this is immediate from
[5, (3.1)]). As $M$ is isomorphic to a direct summand of $\tilde{M}$ our
assumption on $z$ (combined with [5, (3.1)]) implies that $z\in \cal V$.

Let
$$
R^{-}_{2}\, =\, \{-2\epsilon_{2},\, \epsilon_{1}-\epsilon_{2}\}\cup
\{-\epsilon_{2}\pm \epsilon_{j}\, |\, 3\le j\le l\}.
$$
As $p\ne 2$ and $z_{-1}=e_{\gamma}$, there exists $g\in U_1$
such that $\mbox{Supp}\, ((\mbox{Ad}\, g)\cdot z) \cap R^{-}_{2}\, =\,
\emptyset$. Since $\cal V$ is $(\mbox{Ad}\, U_{1})$--invariant we may
(and will) assume that $\mbox{Supp}\, (z)$ does not intersect
with $R^{-}_{2}$.

Let $h_{2}(t)$ be the one-dimensional torus of $L^{\prime}_{1}$
such that
$$
(\mbox{Ad}\, h_{2}(t))\cdot e_{\delta}=
t^{\langle \delta, -2\epsilon_{2} \rangle}\cdot e_{\delta}
$$
for every $\delta\in R$ and every $t\in K^*$. Since
$\langle \delta, -2\epsilon_{2} \rangle \in
\{-2, -1, 0\}$ unless $\delta\in R^{-}_{2}\cup \{\gamma\}$, and
$\langle \gamma, -2\epsilon_{2}\rangle = 1$, we get
$$
t^{2}(\mbox{Ad}\, h_{2}(t))\cdot z = t^{3} e_{\gamma}+
t^{2} y_{2}+t y_{1}+y_{0}
$$
for some $y_{0}, y_{1}, y_{2}\in \frak g$. Hence the Zariski closure of the
conical set
$(\mbox{Ad}\, h_{2}(t))\, (K\, z) \subset\frak g$ contains
$e_{\gamma}$.
As $\cal V$ is conical, Zariski closed and $(\mbox{Ad}\, h_{2}(t))$--stable
we obtain that $e_{\gamma}\in \cal V$.

Let $\beta = -(\epsilon_{1}-\epsilon_{2})$, and let $\frak c$ denote the Lie
subalgebra of $\frak g$ generated by $e_{\beta}$ and $e_{\gamma}$. Since
$p\ne 2$, $\frak c$ is isomorphic to a three-dimensional Heisenberg Lie
algebra. Since $\xi (e_{-\tilde{\alpha}})=\langle e_{\tilde{\alpha}},
e_{-\tilde{\alpha}}\rangle \ne 0$, the derived subalgebra $[{\frak c}, {\frak c}]
\, =\, K\, e_{-\tilde{\alpha}}$ acts invertibly on $\tilde{M}$. Applying
[3, Lemma~2.9] to all composition factors of the $\frak c$--module
$\tilde{M}$ shows that $\tilde{M}$ is free as an $u(e_{\gamma}, \xi)$--module.
By [5, (3.1)], $e_{\gamma}\not\in \cal V$. This contradiction fills the gap
in the proof of the main result of [5] (Theorem~1.1).

We are now going to generalise Theorem~1.1 to the case where the ground field
$K$ has an arbitrary prime characteristic.

\medskip

\noindent
{\bf Theorem~1.1$^{\prime}$}. {\it Let} $G$ {\it be a semisimple and simply
connected algebraic group defined over an algebraically closed field} $K$
{\it of characteristic} $p>0$, {\it and} ${\frak g}\, =\, \mbox{Lie}\, (G)$.
{\it Let} $M$ {\it be a nonzero} $\frak g$--{\it module with}
$p$--{\it character} $\chi\in{\frak g}^*$. {\it Then} $\cal V_{\frak g}(M)
\subseteq\cal N_{p}({\frak g}) \cap {\frak z}_{\frak g}(\chi)$.

\begin{pf}
As in [5] one reduces to the case where $G$ is simple. In view of Theorem~1.1
we may assume that $G\not\cong SL(2)$ and that $p$ is special for $G$.
So either $R$ is of type
${\bf B}_{l}$, ${\bf C}_{l}$, ${\bf F}_4$ and $p=2$ or $R$ is of type
${\bf G}_2$ and $p=3$. Observe that [5, Lemma~2.3(i)] is valid for any $p$
because the $(\mbox{Ad}\, G)$--module $\frak g$ is always isomorphic to
the Weyl module $V(\tilde{\alpha}$). By construction, the bilinear form
$\langle\, \cdot \, ,\, \cdot \, \rangle : {\frak g}\times{\frak g}\rightarrow
K$ is nonzero and $(\mbox{Ad}\, G)$--invariant
(that is [5, Lemma~2.2(i)] holds in our situation).
However, $\mbox{Rad}\, \langle\, \cdot\, ,\, \cdot\, \rangle$ no longer
belongs to the centre of $\frak g$ (that is [5, Lemma~2.2(ii)] fails).

Let $M$ be a nonzero $\frak g$--module with $p$--character $\chi\in {\frak g}^*$
and $z\in \cal V_{\frak g}(M)$. Suppose that $\chi([z, {\frak g}])\ne 0$.
Clearly, this implies $\chi \ne 0$.
First consider the case where $R$ is not of type ${\bf C}_l$. Then
the roots $\alpha_{i_0}$ and $\tilde{\alpha}-\alpha_{i_0}$ are long whence
$[e_{\alpha_{i_{0}}}, e_{\tilde{\alpha}-\alpha_{i_0}}]=\pm e_{\tilde{\alpha}}$.
As $\cal E$ spans $\frak g$ (by [5, Lemma~2.2(i)]) we therefore have
$b_{\chi}(\cal E, \cal E)\ne 0$. Repeating verbatim the argument presented
in [5, PP.~243, 244] one obtains that there exists $e\in \cal E$ such that
$\chi (e)=0$ and $\chi ([z, e])\ne 0$.
Lemma~3.4 of [5] holds in all cases as the proof given in [5, P.~244]
does not require any assumption on $p$. So one can finish the proof
of Theorem~1.1${^\prime}$ for $R\not\cong {\bf C}_l$ as in [5, P.~245].

Finally, suppose $R$ is of type ${\bf C}_l$ and $p=2$. Again ${\frak g} \cong
V(\tilde{\alpha})$ as $(\mbox{Ad}\, G)$--modules whence  $\cal E$ spans
$\frak g$. So there is $e^{\prime}\in \cal E$ with $\chi ([z, e^{\prime}])\ne 0$.
Recall that ${e^{\prime}}^{[p]}=0$. As $p=2$ this yields
$(\mbox {ad}\, e^{\prime})^{2}=0$.
Let ${\frak c}^{\prime}$ denote the Lie subalgebra of $\frak g$ generated by
$e^{\prime}$ and $z$. As $z^{[p]}=0$ one has $(\mbox{ad}\, z)^{2}=0$.
Therefore, ${\frak c}^{\prime}$ is isomorphic to a three-dimensional
Heisenberg Lie algebra. By [5, Lemma~3.4] (which is valid
in all cases), $[z, e^{\prime}]^{[p]}\, =\, \langle z, e^{\prime} \rangle
[z, e^{\prime}]$. It follows that $[z, e^{\prime}]$ acts invertibly on $M$
(one should take into account that $\xi ([z, e^{\prime}])\ne 0$).
Applying [3, Lemma~2.9] to all composition factors of the
${\frak c}^{\prime}$--module $M$ now shows that the $u(z, \chi)$--module $M$
is free. Since this contradicts [5, (3.1)] we are done.
\end{pf}

\medskip

\noindent
{\bf Remark.} In [5], I conjectured that for every $\chi\in {\frak g}^*$
there exists a $\frak g$--module $E$ with $p$--character $\chi$ such that
$\cal V_{\frak g}(E)\, =\, \cal N_{p}({\frak g})\cap {\frak z}_{\frak g}(\chi)$.
This conjecture is proved in [6] under the assumption that the characteristic
of the ground field $K$ is good for $G$. The case of bad characteristics
remains open.

\medskip

{\it Acknowledgements}. I am most thankful to J.~C.~Jantzen, H.~Kraft and
N.~Wallach for their comments on Lemma~3.2 of my paper [5].

\end{document}